\documentclass{amsart}

\title{Generalized String Topology and Derived Koszul Duality}
\author[A.M. Royer]{Aaron M Royer}
\address{Dept. of Mathematics \\
University of Texas \\
1 University Station C1200 \\
Austin, TX 78712-0257}
\email[Royer]
{aroyer@math.utexas.edu}
\date{\today}

\input xy
\xyoption{all}

\usepackage{eucal, bbm}
\usepackage{stmaryrd, hyperref, color, mathrsfs, amssymb, tikz, tikz-cd}
\usetikzlibrary{matrix,arrows,decorations.pathmorphing}

\newtheorem{thm}{Theorem}[section]
\newtheorem{prop}[thm]{Proposition}
\newtheorem{cor}[thm]{Corollary}
\newtheorem{lemma}[thm]{Lemma}

\theoremstyle{definition}
\newtheorem{defn}[thm]{Definition}
\newtheorem{defprop}[thm]{Definition/Proposition}
\newtheorem{const}[thm]{Construction}
\newtheorem{rmk}[thm]{Remark}
\newtheorem{ex}[thm]{Example}

\newcommand{\p}{\partial}
\newcommand{\bbD}{\mathbb{D}}
\newcommand{\bbL}{\mathbb{L}}
\newcommand{\bbP}{\mathbb{P}}
\newcommand{\bbR}{\mathbb{R}}
\newcommand{\bbS}{\mathbb{S}}
\newcommand{\bbU}{\mathbb{U}}
\newcommand{\bbZ}{\mathbb{Z}}
\newcommand{\mtm}{M^{-TM}}

\newcommand{\ntn}{N^{-TN}}
\newcommand{\xtm}{X^{-TM}}
\newcommand{\ytm}{Y^{-TM}}
\newcommand{\xytm}{(X \times_M Y)^{-TM}}
\newcommand{\cA}{\mathcal{A}}
\newcommand{\cC}{\mathcal{C}}
\newcommand{\cD}{\mathcal{D}}

\newcommand{\cS}{\mathcal{S}}
\newcommand{\cT}{\mathcal{T}}
\newcommand{\into}{\hookrightarrow}

\newcommand{\mRk}{m \, \bbR^k}
\newcommand{\nRk}{n \, \bbR^k}
\newcommand{\mnRk}{(m+n) \, \bbR^k}
\newcommand{\id}{\operatorname{id}}
\newcommand{\colim}{\operatorname{colim}}
\newcommand{\hocolim}{\operatorname{hocolim}}
\newcommand{\uhocolim}{\operatorname{uhocolim}}
\newcommand{\sA}{\mathscr{A}}
\newcommand{\sS}{\mathscr{S}}
\newcommand{\sT}{\mathscr{T}}
\newcommand{\sC}{\mathscr{C}}
\newcommand{\sD}{\mathscr{D}}

\newcommand{\sP}{\mathscr{P}}
\newcommand{\sN}{\mathscr{N}}
\newcommand{\Mod}{\operatorname{Mod}}
\newcommand{\Sp}{\operatorname{Sp}}
\newcommand{\Fun}{\operatorname{Fun}}
\newcommand{\FunL}{\operatorname{Fun}^\mathbf{L}}
\newcommand{\Pre}{\operatorname{Pre}}
\newcommand{\End}{\operatorname{End}}
\newcommand{\Aut}{\operatorname{Aut}}
\newcommand{\N}{\mathbf{N}}

\newcommand{\Hom}{\operatorname{Hom}}

\begin{document}

\begin{abstract}
The generalized string topology construction of Gruher and Salvatore assigns to any bundle of $E_n$-algebras $A$ over a closed oriented manifold $M$ a collection of intersection-type operations on the homology of the total space.  These operations are realized by an $H_n$-ring structure on the Thom spectrum $A^{-TM}$ under the Thom isomorphism.  We rigidify and extend this construction to a functor connecting the homotopy theory of spaces and spectra parametrized by $M$ to the homotopy theory of module spectra over the Atiyah-Milnor-Spanier-Whitehead dual $M^{-TM} \simeq \bbD M$.  Then, using an $\infty$-categorical version of Morita theory, we give an alternative description of our construction in terms of the derived Koszul duality (alias bar-cobar duality) between $\Sigma^\infty_+ \Omega M$ and $\bbD M$.
\end{abstract}

\maketitle

\setcounter{tocdepth}{1}
\tableofcontents

\section{Introduction}
\subsection{String Topology}
Let $M^d$ be a closed, connected, smooth manifold.  In \cite{cs99}, Chas and Sullivan introduced a binary operation of degree $-d$ on the homology of $LM = Map(S^1, M)$ called the ``loop product.''  Intuitively, this product mixes the Pontrjagin product on $H_* (\Omega M)$ with the intersection product on $H_* (M)$.  Its original definition of used transversality techniques, and as such was restricted to the case of ordinary homology.  However, the intersection product on $H_* (M)$ may be realized in the stable homotopy category via a Pontrjagin-Thom construction.  The essential point is that a tubular neighborhood $U$ of the diagonal embedding $\Delta \colon M \into M \times M$ is diffeomorphic to the tangent bundle of $M$, and so there is a collapse map
\begin{equation*}
M \times M \to M \times M / ((M \times M) - U) \cong M^{TM}.
\end{equation*}
This collapse induces a commutative and associative product on the Thom spectrum $\mtm$.  To complete the ring spectrum structure, there is a unit $\bbS \to \mtm$ corresponding to the Pontrjagin-Thom collapse familiar from cobordism theory.  With this in mind, Cohen and Jones showed in \cite{cj02} that the loop product is realized by a ring structure on $LM^{-TM}$, the Thom spectrum corresponding to the pullback of the negative tangent bundle of $M$ along evaluation at the basepoint in $S^1$.  An upshot of this construction was that it allowed definition of a loop product in any multiplicative homology theory for which $M$ admits an orientation.  This extended Pontrjagin-Thom construction was further generalized by Gruher and Salvatore in \cite{gs08}, where they showed that for any fiberwise monoid or $E_n$-space $A$ over $M$, the Thom spectrum $A^{-TM}$ is an associative or $E_n$-algebra in the stable homotopy category.  Similarly, given a fiberwise module $B$ over $A$, the Thom spectrum $B^{-TM}$ is a module spectrum over $A^{-TM}$.  However, the target of their constructions is the classical stable homotopy category, where it is known that ring and module theory behave less than ideally.  One of our main results is a lift of their work to a modern symmetric-monoidal point-set category of spectra, and with this extra structure we are able to connect with other parts of parametrized homotopy theory.

\subsection{Rigidifying the Pontrjagin-Thom Construction}
The first step toward lifting products on the Thom spectra $A^{-TM}$ to a homotopy-coherent setting is rigidifying the product on $\mtm$.  This problem was attacked successfully by Cohen in \cite{c04}, where he realized $\mtm$ as a commutative symmetric ring spectrum at the cost of building in coherence machinery to track embedding and tubular neighborhood data.  In fact, there are two models constructed in that work.  One has nice functoriality properties resulting from independence from embedding data, but is not unital in any reasonable sense.  The other is significantly smaller and has a strict unit given by the Pontrjagin-Thom collapse, but is not really functorial for maps more general than embeddings of closed submanifolds.  Our first result is a simplification of Cohen's second construction, and is the base case for the point-set portion of our work.

\begin{thm}
There is a commutative symmetric ring spectrum $\mtm$ and a morphism of symmetric spectra
\begin{equation*}
\mtm \wedge M_+ \to \bbS
\end{equation*}
inducing an equivalence of commutative symmetric ring spectra $\mtm \simeq F(M_+, \bbS)$.  Moreover, this agrees with the equivalence of \cite{a61} up to homotopy.
\end{thm}

In fact, there are many such symmetric ring spectra, as the construction depends on an embedding of $M$ and a tubular neighborhood of that embedding.  From our perspective, a main advantage of the above theorem is that it shows none of these choices matter up to equivalence of $E_\infty$-rings.  Just as in Cohen's case, however, this construction has fairly limited functoriality in $M$.  For our immediate purposes this will not matter too much, since we are mainly concerned with spaces parametrized by a single, fixed manifold.  This compromise allows for small models of the resulting Thom spectra which are relatively easy to work with.

\subsection{Generalized String Topology}
Next, we use the techniques involved in our construction of $\mtm$ as a commutative symmetric ring spectrum to produce a functor from the category $\sT_{/M}$ of spaces over $M$ to modules over $\mtm$ which rigidifies the previous constructions of Thom spectra $A^{-TM}$.  We call this functor generalized string topology, and abusively denote it by $(-)^{-TM}$.  Like the Atiyah duality construction above, $(-)^{-TM}$ depends on an embedding of $M$ into a Euclidean space and a tubular neighborhood of that embedding.  Once these choices are made, our functor enjoys many nice properties.  We are able to show the following.

\begin{thm}
The generalized string topology functor
\begin{equation*}
(-)^{-TM} \colon \sT_{/M} \to \Sigma \sS_{\mtm}
\end{equation*}
is a topological lax symmetric monoidal left adjoint and preserves weak equivalences and tensors with spaces.
\end{thm}

Unfortunately, the construction does not interact well with cellular constructions, and appears not to be a left Quillen adjoint.  For this reason, it is inconvenient to use classical model category theory to analyze the generalized string topology functor. Luckily, in the past few years, new and more flexible categorical foundations for homotopy theory have appeared.  In particular, we will make use of the theory of $\infty$-categories, and by translating into this setting we obtain a simple conceptual characterization of $(-)^{-TM}$ independent of any choices.

\subsection{$\infty$-categorical Context}

The theory of $\infty$-categories (or, more accurately, $(\infty,1)$-categories) may be thought of as category theory ``up to coherent homotopy.''  Roughly, an $\infty$-category $\cC$ consists of a collection of objects, and for any pair of objects $X,Y$ a space $\cC(X,Y)$ of morphisms from $X$ to $Y$.  There is a composition operation defined up to homotopy which is associative and unital up to all higher coherence homotopies.  Conceptually, $\infty$-category theory is a simultaneous generalization of ordinary category theory and homotopy theory of spaces, and in fact strictly contains both.  Thanks to the work of many authors over many years, we have several equivalent concrete models of $\infty$-categories and a large tool kit for manipulating them.  The resulting theory formally looks very similar to ordinary category theory.  For example, there are good concepts of limits, colimits, adjoint functors, (co)monads, (symmetric) monoidal structures, etc.

A large, important class of $\infty$-categories arise as intermediate between model categories, (or more generally Waldhausen or homotopical categories,) and their homotopy categories.  Passing from a model category to its ``underlying'' $\infty$-category forces weakly equivalent objects to be isomorphic while retaining higher homotopical coherences.  Moreover, this passage sometimes allows for analysis of functors between model categories which are not Quillen but nevertheless have some nice homotopical properties.  A primary advantage from our point of view is that functors which pass to underlying $\infty$-categories and preserve homotopy colimits induce functors that preserve colimits on the $\infty$-categorical level.  With this in hand, Theorem 1.2 now takes the following form.

\begin{thm}
The generalized string topology functor
\begin{equation*}
(-)^{-TM} \colon \cT_{/M} \to \Mod_{\mtm}
\end{equation*}
is lax symmetric monoidal and preserves colimits.
\end{thm}

The change in notation is to emphasize that we are now dealing with $\infty$-categories.  Proving the colimit statement is perhaps the most technically demanding part of our work, and involves the topological enrichment in a critical way.

There are two classes of $\infty$-categories whose properties play an essential role in the formulation and proof of our final result.  The first class consists of the presheaf $\infty$-categories of the form $\Pre(\cC) := \Fun(\cC^{op}, \cT)$ for some small $\infty$-category $\cC$.  Just as is the case in ordinary category theory, there is a fully-faithful Yoneda embedding $\Pre(\cC)$ given by sending an object to the functor it represents.  The category $\Pre(\cC)$ is in a precise sense the free cocompletion of $\cC$.  This means, for example, that for any $\cD$ the Yoneda embedding determines an equivalence of $\infty$-categories
\begin{equation*}
\FunL(\Pre(\cC), \cD) \simeq \Fun(\cC, \cD).
\end{equation*}
Here, $\FunL$ denotes the $\infty$-category of colimit-preserving functors.  As a concrete example, suppose $\cC$ has a single object whose endomorphisms are a topological group $G$ (or more generally a grouplike $A_\infty$-space $\Omega X$.)  Then $\Pre(\cC)$ is the $\infty$-category $\cT^{\, G}$ of spaces with a $G$-action, considered with the Borel homotopy theory (i.e. equivalences are created through the forgetful functor to spaces.)

The other important class of $\infty$-categories for our purposes consists of the stable $\infty$-categories.  These are $\infty$-categories with properties abstracted from the $\infty$-category of spectra, over which they are all canonically enriched.  The role they play in the general $\infty$-category theory is analogous to that of abelian categories in the ordinary theory, though they more closely resemble triangulated categories.  More precisely, the homotopy category of a stable $\infty$-category is triangulated, and considering a stable $\infty$-category rather than its homotopy category frequently leads to better technical behavior.  A typical example is the $\infty$-category $\Mod_R$ of right modules over an $A_\infty$-ring $R$, or more generally the $\infty$-category of presheaves of spectra on a small $\infty$-category $\cC$.  Under mild technical hypotheses, one may ``freely stabilize'' an $\infty$-category $\cC$ by forming the $\infty$-category of spectrum objects in $\cC$, denoted $\Sp(\cC)$.  This comes equipped with an adjunction
\begin{equation*}
\Sigma_\cC^\infty \colon \cC \rightleftarrows \Sp(\cC) \colon \Omega^\infty_\cC
\end{equation*}
with the property that precomposition with $\Sigma_\cC^\infty$ induces an equivalence of $\infty$-categories
\begin{equation*}
\FunL(\Sp(\cC), \cA) \simeq \FunL (\cC, \cA)
\end{equation*}
for any stable $\infty$-category $\cA$.  In the case of $\cC = \cT{\, G}$ from above, we may identify $\Sp(\cC)$ with the $\infty$-category $\Mod_{\Sigma^\infty_+ G}$ of modules for the spherical group ring of $G$.

\subsection{Parametrized Homotopy Theory, Derived Morita Theory and the Comparison Theorem}

To arrive at our main result, we use $\infty$-categorical analogues of two classical results.  The first is an alternative characterization of the homotopy theory of parametrized spaces originally due to Dror, Dwyer and Kan.  Suppose that $X$ is a connected, based CW complex.  Then it is classical that the category of local coefficient systems with values in an abelian category $\sA$ is equivalent to the category of representations of $\pi_1 (X)$ into $\sA$, the equivalence being given in one direction by taking the fiber at the basepoint and in the other by an associated bundle construction.  This generalizes very nicely once one allows the full based loop space $\Omega X$ to act rather than just its group of path components.  In particular, we shall use the following special case.

\begin{thm}
Let $X$ be a connected, based CW complex.  Then there is an equivalence of $\infty$-categories
\begin{equation*}
(-) \times_X * \colon \cT_{/X} \rightleftarrows \cT^{\, \Omega X} \colon (-) \times_{\Omega X} *
\end{equation*}
Here the left-hand side is the $\infty$-category of spaces parametrized over $X$, and the right-hand side is the $\infty$-category of spaces equipped with an action of the group-like $A_\infty$-space $\Omega X$.
\end{thm}

The notation above may be slightly misleading, but is reasonable from the $\infty$-categorical standpoint.  Each functor is actually the derived functor of what is written on the point-set level.  To obtain a $\Omega X$-space from a space over $X$, for instance, one takes the fiber product with the based path space of $X$, a fibrant replacement of $* \to X$.  A Quillen equivalence of model categories presenting this equivalence, as well as its relation to the theory of homotopy sheaves, may be found in \cite{shu08}.  Note that there is an equivalence between the $\infty$-category of spaces over $X$ and comodules for $X$, considered as a commutative coalgebra under the diagonal map.  With this in mind, the above equivalence may be thought of as an unstable or nonabelian bar-cobar duality.

The other classical result to generalize is the Eilenberg-Watts theorem.  In its original form, the theorem states that for rings $R,S$ any colimit-preserving functor $\Mod_R \to \Mod_S$ is naturally equivalent to $(-) \otimes_R P$ for an $R-S$-bimodule $P$.  This theorem and even its proof carry over essentially verbatim to the $\infty$-categorical setting.

\begin{thm}
Let $R, S$ be $A_\infty$-rings, $F \colon \Mod_R \to \Mod_S$ a colimit-preserving functor.  Then there is an $R-S$-bimodule $P$ and a natural equivalence
\begin{equation*}
F \simeq (-) \wedge_R P
\end{equation*}
\end{thm}

The proof is actually just a special case of the discussion at the end of the previous section, as $\Mod_R$ may be identified with an $\infty$-category of presheaves of spectra.  Putting it all together, we obtain the following ``algebraic'' construction of generalized string topology, our final result.  We write $S^{-TM}$ for $S^{-d}$ with the $\Sigma_+^\infty\Omega M$-module structure determined by the classifying map $M \simeq B \Omega M \to BhAut(\bbS)$ classifying the stable spherical fibration corresponding to the stable normal bundle $-TM$.
\begin{thm}
Under the above equivalence $\cT_{/M} \simeq \cT^{\, \Omega M}$, the generalized string topology functor $(-)^{-TM}$ is naturally equivalent to the functor
\begin{equation*}
(-) \wedge_{\Sigma_+^\infty \Omega M} S^{-TM} \colon \Mod_{\Sigma^\infty_+ \Omega M} \to \operatorname{Mod}_{\bbD M}
\end{equation*}
restricted to suspension spectra.
\end{thm}

The proof amounts to computing the value of each functor on the image of the Yoneda embedding.  This amounts to $\Omega M$ acting on itself on the one side, and $* \simeq PM \to M$ on the other.  We derive the title of our work from the fact that $\bbD M$ is the Koszul dual of $\Sigma_+^\infty \Omega M$, i.e.
\begin{equation*}
\bbD M \simeq \End_{\Sigma^\infty_+ \Omega M} (\bbS).
\end{equation*}
At this point we note that we obtain model-independence up to essentially unique equivalence as a consequence of this result.  By contrast, a direct comparison between functors based on different embeddings and tubular neighborhoods is significantly more intricate.

\subsection{Further Directions}

A significant portion of the work following Chas and Sullivan's original preprint has been on topological quantum field-theoretic aspects of string topology.  It is shown there that $H_{*-d} (LM;k)$ carries a Batalin-Vilkovisky algebra structure, built out of the loop product and the $S^1$-action.  In fact, it can be shown that this structure arises (upon the passage to homology) from a topological conformal field theory.  From \cite{cos07}, one supposes that this TCFT is the closed-string portion of an open-closed theory associated to an $A_\infty$-category of paths in the manifold.  A construction of a category of this kind was sketched in \cite{bct09}, and its relation to symplectic topology of cotangent bundles explored.  The coherence machinery developed in this work can be used to give a precise construction of the string topology category of \cite[2.12]{bct09}.  We intend to give this construction and explore its consequences in a future paper.

Another feature of our constructions here is the formal extension of string topology operations to parametrized spectra, rather than just spaces.  The theory developed in \cite{ms06} and \cite{abghr} shows that this is the correct conceptual setting for twisted (co)homology, Thom spectra and orientations.  Therefore, our technology extends to give string topology operations on twisted generalized homology of spaces parametrized by closed smooth manifolds.  We will give details in a sequel.

Finally, in the case that $M$ comes equipped with the action of a finite group $G$, the arguments of this work go through essentially unchanged to obtain a generalized string topology functor from $G$-spaces over $M$ to a category of genuine $G$-equivariant symmetric spectra, such as those constructed in \cite{man04}.  The first appendix of that work also contains a version of genuine $G$-equivariant symmetric spectra for a general compact Lie group $G$, but it is somewhat more subtle than the finite group case.  In a sequel, we will give details to extend all of the constructions of this work to the general equivariant case.

\subsection{Notation}

As detailed above, we will use a combination of model categories and $\infty$-categories in this work.  To avoid confusion, we adopt the notational convention that ordinary categories (including model categories) are nearly always denoted with script letters such as $\sC$, $\sS$ and $\sT$, and calligraphic letters such as $\cC$, $\cS$ and $\cT$ are reserved for $\infty$-categories.  We will write $\bbS$ for the sphere spectrum independent of the categorical context, but decorate it to refer to specific models.

\subsection{Acknowledgements}
The author would like to thank David Ayala, Ralph Cohen, Matt Pancia,  Hiro Tanaka, and Sam Taylor for helpful conversations at various stages of this project.  He is also particularly grateful to Andrew Blumberg for all the support and patience.

\section{Preliminaries}
\subsection{Model Categories}
The domain category for our generalized string topology functor is the over-category $\sT_{/M}$ of (compactly generated) spaces mapping to $M$, a fixed base manifold.  This inherits a model category structure from any model structure on the underlying category $\sT$ of spaces (via the forgetful functor.)  We will work in the following model category of parametrized spaces.

\begin{defprop}
Define the $m$-model structure on $\sT_{/M}$ as the model structure in which:
\begin{itemize}
\item Weak equivalences are those maps over $M$ which are weak homotopy equivalences on total spaces,
\item Fibrations are Hurewicz fibrations over $M$, and
\item The cofibrant objects in are spaces over $M$ with the homotopy type of CW complexes.
\end{itemize}
\end{defprop}

The $m$ stands for ``mixed,'' and was originally developed by Michael Cole.  See \cite[\S 3]{shu08} for an exposition of the $m$-model structure, and \cite[Part II]{ms06} for complete details.  We will consider $\sT_{/M}$ with the Cartesian symmetric monoidal structure.  The reader should be warned that with this product $\sT_{/M}$ is not a monoidal model category in the usual sense, but preserves all weak equivalences as long as one argument is fibrant.  As we will eventually be concerned with its underlying $\infty$-category, we merely remark that restriction to cofibrant-fibrant objects will suffice for our purposes.  The issue will be dealt with more concretely below.  We also note that this monoidal structure is not closed because of pathologies in point-set topology.  This need not concern us in this work, but see the sources cited above for full details.

The target of our functor is a symmetric monoidal category of spectra where structure maps are defined for $k$-fold suspensions, typically for $k \gg 1$.  We construct such categories in the appendix and compare them to the ordinary categories of diagram spectra from \cite{mmss}.  We will state the basic definitions and theorems we require, and defer proofs to the appendix.

\begin{defn}
Let $\Sigma$ be the symmetric monoidal groupoid with object set $\bbZ_{\geq 0}$, morphisms defined by
\begin{equation*}
\Sigma(m,n) = \begin{cases} \Sigma_m & m=n \\ \emptyset & m \neq n \end{cases},
\end{equation*}
and symmetric monoidal structure addition on objects and block-sum on morphisms.  We implicitly consider $\Sigma$ as a discrete topological category.  Define a symmetric sequence (of spaces) to be a (topological) functor $\Sigma$ to (compactly generated) topological spaces, and denote the category of such functors and natural transformations $\Sigma \sT$.  Similarly define the category of symmetric sequences of based spaces, denoted $\Sigma \sT_*$.
\end{defn}

By-now standard categorical machinery yields the following.

\begin{prop}
The categories $\Sigma \sT$ and $\Sigma \sT_*$ are closed symmetric monoidal under Day convolution.
\end{prop}

Write $\times$ for the monoidal product in $\Sigma \sT$, and $\wedge$ for the product in $\Sigma \sT_*$.  We refer to \cite[\S21]{mmss} for complete details of their construction, but list the universal property of $\wedge$ here as we will use it repeatedly.

\begin{prop}
Let $X, Y, Z \in \Sigma \sT_*$.  Then $Map (X \wedge Y, Z)$ is naturally isomorphic to the space of collections of maps $\{X(m) \wedge Y(n) \to Z(m+n)\}_{m,n \in \bbZ_{\geq 0}}$ which are $\Sigma_m \times \Sigma_n$-equivariant.
\end{prop}
There is a similar universal property for the product of unbased symmetric sequences.  These serve the purpose of transporting associativity coherence questions (e.g. Mac Lane's pentagon) to the corresponding structures in $\sT_*$ and $\sT$, respectively, where they are easily dealt with by standard techniques.  The symmetry structure and coherences are more subtle.  To understand them, it is useful to write an explicit formula for the spaces in a smash product of symmetric sequences.  To this end, for $X, Y \in \Sigma \sT_*$ we have
\begin{equation*}
(X \wedge Y)(r) \cong \bigvee_{m+n=r} (\Sigma_r)_+ \wedge_{\Sigma_m \times \Sigma_n} X(m) \wedge Y(n).
\end{equation*}
On wedge factors, the twist isomorphism $\tau_{X,Y} \colon X \wedge Y \to Y \wedge X$ is the twist isomorphism
\begin{equation*}
\tau_{X(m), Y(n)} \colon X(m) \wedge Y(n) \to Y(n) \wedge X(m)
\end{equation*}
of based spaces on the right-hand smash factors but acts on $\Sigma_r$ by multiplication by the block transposition $\chi_{m,n} := (m+1, ..., n, 1, ..., m)$.  As a result, to give a commutative monoid in symmetric sequences using Proposition 2.3 requires an additional $\Sigma_2$-equivariance, and similar data must be given on functors and natural transformations for them to respect the symmetric monoidal structure.

For each $k \geq 1$, we define a commutative monoid $\bbS_k \in \Sigma \sT_*$ by
\begin{equation*}
\bbS_k (m) = S^k \wedge ... \wedge S^k \cong S^{mk},
\end{equation*}
where $\Sigma_m$ acts by permutation of the $S^k$ factors.  Here we choose a concrete model for $S^k$ as the one-point compactification of $\bbR^k$, and with this a homeomorphism of $\bbS_k (m)$ with the one-point compactification of $\bbR^m \otimes \bbR^k \cong \bbR^{mk}$.  Under this identification, the symmetric group $\Sigma_m$ acts by block permutation matrices, and in particular by linear isometries.  Write $\Sigma \sS_k$ for the symmetric monoidal category of modules over $\bbS_k$, which we will often call $k$-symmetric spectra.  Note that when $k=1$ this is the usual category $\Sigma \sS$ of symmetric spectra of topological spaces.  Using the technology of \cite[\S2 \& 6]{mmss}, we may endow the category $\Sigma \sS_k$ with a compactly generated topological model structure.

\begin{defprop}
Define the level model structure on $\Sigma \sS_k$ as the model structure in which:
\begin{itemize}
\item Weak equivalences are those morphisms of $k$-symmetric spectra which are weak homotopy equivalences at each space in the sequence,
\item Fibrations are morphisms which are Serre fibrations at each space in the sequence,
\item Cofibrations are morphisms which have the left-lifting property with respect to morphisms which are acyclic Serre fibrations at each space.
\end{itemize}
\end{defprop}

The stable model structure is obtained as a left Bousfield localization of this model structure.  It is determined up to Quillen equivalence by requiring that the fibrant objects are precisely the collection $C$ of $\Omega^k$-spectra, which we define to be those $k$-symmetric spectra $E$ for which the adjoint structure maps
\begin{equation*}
E(m) \to \Omega^k E(m+1)
\end{equation*}
are weak homotopy equivalences for all $m$.  The weak equivalences in this stable model structure are then morphisms of $k$-symmetric spectra which are $C$-local equivalences.  We will often make use of a stricter version of weak equivalence.
\begin{defn}
Let $E \in \Sigma \sS_k$.  Then the homotopy groups of $E$ are defined by
\begin{equation*}
\pi_s (E) := \colim_{r \to \infty} \pi_{rk + s} E(r), \ s \in \bbZ.
\end{equation*}
Here, the colimit is through the structure maps of the $k$-symmetric spectrum $E$.  A $\pi_*$-isomorphism of $k$-symmetric spectra is a morphism $f \colon E \to E'$ which induces isomorphisms $\pi_s (E) \cong \pi_s (E')$ for all $s \in Z$.
\end{defn}
The homotopy groups of $k$-symmetric spectra suffer from the same strangeness with respect to stable equivalences as in the case $k=1$.  Suffice it to say that the $\pi_*$-isomorphisms are contained in the weak equivalences of the stable model structure on $\Sigma \sS_k$.  The most important result for our purposes is the following comparison theorem.

\begin{prop}
For each $k$, the forgetful functor from symmetric spectra to $k$-symmetric spectra is the right adjoint in a Quillen equivalence
\begin{equation*}
\bbP_k \colon \Sigma \sS \rightleftarrows \Sigma \sS_k \colon \bbU_k
\end{equation*}
The left adjoint $\bbP_k$ is strong symmetric monoidal, and $\bbU_k$ is lax symmetric monoidal.  In particular, $\Sigma \sS_k$ is a symmetric monoidal model of the stable homotopy category.
\end{prop}

Now suppose that $A$ is a commutative monoid in $\Sigma \sS_k$, or equivalently a commutative $\bbS_k$-algebra in $\Sigma \sT_*$.  Then we may define a symmetric monoidal category of $k$-symmetric module spectra over $A$, denoted $\Sigma \sS_A$.  There is an obvious forgetful functor $\bbU \colon \Sigma \sS_A \to \Sigma \sS_k$.  In parallel to the classical case of commutative rings, $\bbU$ has a left adjoint free $A$-module functor.  We use the right adjoint $\bbU$ to create a stable model structure on $\Sigma \sS_A$.  That is, fibrations and weak equivalences of $A$-module spectra coincide with fibrations and weak equivalences of underlying $k$-symmetric spectra, and cofibrations are those $A$-module morphisms with the appropriate lifting properties.  With this model structure, $\bbU$ is obviously right Quillen, and so together with the free functor we have a Quillen adjunction.  More generally, we have the following.

\begin{prop}
Let $f \colon B \to A$ be a morphism of $\bbS_k$-algebras.  Then restricting and extension of scalars form a Quillen adjunction
\begin{equation*}
f_* \colon \Sigma \sS_B \rightleftarrows \Sigma \sS_A \colon f^*
\end{equation*}
If $f$ is an equivalence, then this is a Quillen equivalence.  If $A,B$ are commutative $\bbS_k$-algebras, then $f_*$ is strong symmetric monoidal and $f^*$ is lax symmetric monoidal.
\end{prop}

As we discuss in the Appendix, there is nothing special to $\Sigma$ in these constructions, and any of the standard diagrams will work.  In fact, our proof of Proposition 2.3 actually involves $k$-orthogonal spectra.  However, it turns out our later work is rather specific to the case of symmetric spectra, as we will see in Section 3.

\subsection{$\infty$-Categories}

While many of our technical constructions take place in specific model categories of spaces and spectra, the conceptual heart of this work is best expressed in the setting of $\infty$-categories (or, more precisely, $(\infty,1)$-categories.)  Working in this theory, all functors are natively derived, $\infty$-categorical limits and colimits correspond to point-set homotopy limits and colimits, uniqueness is expressed by the phrase ``up to contractible choice,'' composition of morphisms is associative up to coherent homotopies, et cetera.  These properties all contribute to the formation of a good setting for making formal arguments in modern stable homotopy theory, as well as other ``higher'' or ``derived'' geometric settings.  There are by now several excellent models for the homotopy theory of $\infty$-categories, and the work of several authors have shown them all to be equivalent.  For the sake of concreteness, we will use the quasicategories of \cite{jo02}, \cite{htt}, \cite{ha}, as their theory is very well-developed.

There are several methods of obtaining quasicategories from model categories, with different functoriality properties.  The most basic is the following.  Suppose given a model category $\sC$ with subcategory of weak equivalences $W_\sC$.

\begin{defn}
The localization $\N \sC [W_\sC^{-1}]$ is the initial quasicategory receiving a functor from $\N \sC$ with the property that every arrow in $W$ is sent to an equivalence.
\end{defn}

This construction is conceptually quite satisfying, but in practice it suffers from poor functoriality properties.  In particular, for a functor between model categories $F \colon \sC \to \sD$ to induce a functor $F \colon \N \sC[W_\sC^{-1}] \to \N \sD[W_\sD^{-1}]$, one needs preservation of all weak equivalences.  This excludes a great many functors of interest.  Luckily, one does not need the full model category to present its underlying $\infty$-category.  For example, one can restrict to full subcategories of cofibrant or fibrant objects to present equivalent $\infty$-categories and achieve good behavior with respect to Quillen functors, respectively.  More generally, replacing $\sC$ with any deformation retract in the sense of \cite[39.2]{dhks} will yield an equivalent $\infty$-category, and therefore deformable functors induce functors between underlying $\infty$-categories. (See \cite[VII]{dhks} for a thorough development of the theory of homotopical categories and deformable functors, and \cite[1.3.4.1]{ha} for details on the localization procedure.)

In nature, one often works with model categories equipped with a compatible simplicial or topological enrichment.  The theory of simplicial categories is known to be a model of the theory of $\infty$-categories.  More precisely, the model categories of simplicial categories and of quasicategories are Quillen equivalent via an enriched version of the nerve construction.  This gives another method of obtaining the $\infty$-category underlying a simplicial model category, after restricting to the full simplicial subcategory consisting of cofibrant-fibrant objects.  Fortunately, this method is equivalent to the localization construction above by \cite[1.3.4.20]{ha}.  This method will be important in Section 6, when we use topological enrichment to show that generalized string topology preserves homotopy colimits.

The theory of symmetric monoidal $\infty$-categories is somewhat more intricate than that of ordinary symmetric monoidal categories.  Roughly, allowing for higher homotopies leads to the need for an infinite hierarchy of coherences beyond Mac Lane's pentagon, and the resulting objects are in fact a generalization of infinite loop spaces (via the equivalence between $\infty$-groupoids and homotopy types.)  A good theory has been developed, and through the methods mentioned above one may obtain symmetric monoidal $\infty$-categories from symmetric monoidal model categories \cite[\S 4.1.3]{ha}.  Some care must be taken, since the monoidal product in most examples does not preserve all weak equivalences in each variable separately, and therefore must be derived.  However, in the case of a monoidal model category, restricting to cofibrant objects is sufficient.  The derived monoidal product $\otimes^\bbL$ is then a left derived functor.  Therefore, in the presence of functorial factorizations, the derived monoidal product comes with a natural transformation to the underived monoidal product, and so lax (symmetric) monoidal functors that preserve weak equivalences between cofibrant objects induce lax (symmetric) monoidal functors on the underlying $\infty$-categories.  Additionally, (commutative) monoids in $\sC$ pass to (commutative) monoids in $\cC$, whether or not they are cofibrant.

In the latter portion of our work, we will encounter $\infty$-categories of modules over ring spectra.  These have additional structure reflecting their ``algebraic'' origin.  Most notably, they are stable in the following sense.

\begin{defn}
A pointed $\infty$-category $\cC$ is an $\infty$-category with a distinguished zero object.  A stable $\infty$-category $\cA$ is a pointed $\infty$-category with all finite limits and colimits and such that a square diagram $\Delta^1 \times \Delta^1 \to \cA$ is Cartesian if and only if it is coCartesian.
\end{defn}

It follows from these axioms that the homotopy category of a stable $\infty$-category is triangulated, with triangles coming from (co)Cartesian squares where one of the ``middle'' corners is the zero object.  If both middle corners are the zero object, then the square represents the suspension/shift functor, which is seen to be an equivalence.  It is a nice exercise to show that stable model categories (in the sense of e.g. \cite{ss03}) have stable underlying $\infty$-categories though one of the constructions above.  Given a general $\infty$-category $\cC$ satisfying mild hypotheses, (namely presentability,) one may construct a stable $\infty$-category $\Sp(\cC)$ of spectrum objects in $\cC$ by mimicking the construction of the category of spectra out of the category of spaces.  One obtains an adjunction
\begin{equation*}
\Sigma^\infty_+ \colon \cC \rightleftarrows \Sp(\cC) \colon \Omega^\infty.
\end{equation*}
This adjunction manifests a universal property for $\Sp(\cC)$ in the following sense.  For $\infty$-categories $\cC$, $\cD$, write $\FunL(\cC,\cD)$ for the full sub-$\infty$-category of functors spanned by those that preserve colimits.

\begin{prop}[{\cite[1.4.4.5]{ha}}]
Suppose $\cC$ is presentable and $\cD$ is both presentable and stable.  There is an equivalence of $\infty$-categories
\begin{equation*}
\FunL(\Sp(\cC), \cD) \simeq \FunL(\cC, \cD)
\end{equation*}
induced by pre-composition with $\Sigma_+^\infty \colon \cC \to \Sp(\cC)$.
\end{prop}

That is, $\Sigma^\infty_+$ is the initial colimit-preserving functor from $\cC$ to a stable $\infty$-category.  We will exploit this property to extend the generalized string topology functor from parametrized spaces to parametrized spectra, and explore the consequences of this extension in a sequel.  Details of the construction of $\Sp(\cC)$ may be found in \cite[\S 1.4.2]{ha}.

To fix notation for later, we write $\cT$ for the $\infty$-category of spaces, and $\cS$ for the $\infty$-category of spectra.  We call monoids in $\cS$ $A_\infty$-rings and commutative monoids $E_\infty$-rings, as the $\infty$-categories of these objects are equivalent to the $\infty$-categories underlying point-set models for $A_\infty$ and $E_\infty$ ring spectra, respectively.  For an $A_\infty$-ring $A$ we write $\operatorname{LMod}_A$ for the $\infty$-category of left $A$-modules.  If $A$ is an $E_\infty$-ring, then we write $\operatorname{Mod}_A$ for its symmetric monoidal $\infty$-category of left (or equivalently right) modules.  These module categories are constructed from scratch in \cite[\S 4.4]{ha}, or can be obtained from suitable model categories of modules, such as the ones constructed in \cite[\S 12]{mmss}.

\section{Multiplicative Atiyah Duality}

Let $M$ be a closed connected smooth manifold.  In this section, we give a model for the Thom spectrum of the stable normal bundle of $M$ as a commutative $k$-symmetric ring spectrum, and compare it to the Spanier-Whitehead dual $F(M, \bbS_k)$.

Fix an embedding $e \colon M \into \bbR^k$, and write $\mRk$ for the tensor product $\bbR^m \otimes \bbR^k$.  The space $\mRk$ carries an action of the symmetric group $\Sigma_m$ by permutation of the standard basis vectors.  Denote the diagonal linear map $\bbR^k \into \mRk$ by $\Delta_m$, and observe that it is a conformal embedding with scaling factor $\sqrt{m}$.  We will write $e_m$ for the composition $\Delta_m \circ e$.  By the tubular neighborhood theorem, there exists $\epsilon > 0$ with the following two properties:

\begin{enumerate}
\item $\overline{B}_\epsilon (e(M))$ is diffeomorphic to a unit disk bundle in the normal bundle $\eta_e$.
\item $\overline{B}_\epsilon (e_2(M)) \cap e(M) \times e(M)$ is a tubular neighborhood of the diagonal $M \into M \times M$, and hence is diffeomorphic to a unit disk bundle in the tangent bundle $TM$.
\end{enumerate}

Condition (2) implies the appropriate analogue for any iterated diagonal embedding.  Let $L_e$ be the maximum such number, and fix $\epsilon \in (0, L_e)$.  Write $D (m)$ for $\overline{B}_\epsilon (e_m(M))$, and observe that it is a $\Sigma_m$-invariant submanifold-with-boundary of $\mRk$.  Write $S (m)$ for $\p D(m)$.  Note that if we allow $m$ to vary through the natural numbers, $D(-)$ and $S(-)$ define symmetric sequences.  They will be important in Section 5, and referred to by $D$ and $S$, respectively.

\begin{defn}
The symmetric sequence $\mtm$ is defined by
\begin{align*}
\mtm (0) &= S^0 \\
\mtm (m) &= D (m) / S (m), \ m \geq 1,
\end{align*}
with $\Sigma_m$-action inherited from $D(m)$.
\end{defn}

Remarkably, this symmetric sequence is a highly structured model for the Atiyah dual of $M$, as we now show.

\begin{prop}
The symmetric sequence $\mtm$ is a commutative $\bbS_k$-algebra via the Pontrjagin-Thom construction.
\begin{proof}
We first define a commuative multiplication
\begin{equation*}
\mu \colon \mtm \wedge \mtm \to \mtm.
\end{equation*}
For now, $\wedge$ means smash product of symmetric sequences.  By Proposition 2.2, this is equivalent data to the specificiation of $\Sigma_m \times \Sigma_n$-equivariant maps
\begin{equation*}
\mu \colon \mtm (m) \wedge \mtm (n) \to \mtm (m+n)
\end{equation*}
satisfying appropriate associativity and commutativity constraints.  For $m = 0$ or $n = 0$, use the unit isomorphism $X \wedge S^0 \cong X \cong S^0 \wedge X$.  Otherwise, note that the left-hand space is homeomorphic to
\begin{equation*}
(D (m) \times D (n)) / (S (m) \times D (n)) \cup (D (m) \times S (n)),
\end{equation*}
and the right-hand space is defined to be
\begin{equation*}
D(m+n) / S(m+n).
\end{equation*}
We have inclusions
\begin{equation*}
D(m+n) \subset D(m) \times D(n)
\end{equation*}
and
\begin{equation*}
D(m+n) \cap (S (m) \times D (n)) \cup (D (m) \times S (n)) \subset S(m+n).
\end{equation*}
Thus the quotient map $D(m+n) \to \mtm(m+n)$ extends to all of $D(m) \times D(n)$, and this extension factors through the quotient map for $\mtm (m) \wedge \mtm (n)$.  Diagrammatically,
\begin{equation*}
\begin{tikzcd}
D(m) \times D(n) \arrow{d} \arrow[dashed]{dr} \arrow[hookleftarrow]{r} & D(m+n) \arrow{d} \\
\mtm(m) \wedge \mtm(n) \arrow[dashed]{r} & \mtm(m+n)
\end{tikzcd}
\end{equation*}
The bottom arrow is our multiplication map $\mu$.  Every map in this diagram is $\Sigma_m \times \Sigma_n$-equivariant.  Associativity of iterated applications of $\mu$ follows from the universal property of quotients.  Commutativity follows from the fact that every map in the above diagram is also $\Sigma_2$-equivariant, where the nontrivial element in $\Sigma_2$ acts on the left spaces by permuting factors and the right spaces by the map exchanging $\mRk$ and $\nRk$ in $\mnRk$.

Next, we give a multiplicative map $\eta \colon \bbS_k \to \mtm$ which will act as the unit and endow $\mtm$ with an $\bbS_k$-algebra structure.  Recall that we have chosen a concrete model of $\bbS_k (m)$ as the one-point compactification of $\mRk$.  For each $m \geq 1$, let $\bbS_k (m) \to \mtm (m)$ be the Thom collapse associated to $D(m) \subset \mRk$.  This collapse is $\Sigma_m$-equivariant since $D(m)$ is an invariant subspace of $\mRk$.  Compatibility of this map with multiplication is a consequence of the fact that one-point compactification sends Cartesian products to smash products, together with the universal property of quotients.  From this description, we see as a consequence of associativity of $\mu$ that the two action maps
\begin{equation*}
\mtm \wedge \bbS_k \wedge \mtm \rightrightarrows \mtm \wedge \mtm
\end{equation*}
agree after multiplication.  This means that $\mu$ factors through the balanced smash product $\mtm \wedge_{\bbS_k} \mtm$.  A check of the definitions shows that the action map $\bbS_k \wedge \mtm \to \mtm$ is compatible with the twist isomorphism, and therefore $\mtm$ has the structure of a commutative $\bbS_k$-algebra.
\end{proof}
\end{prop}

\begin{rmk}
We can see at this stage why $k$-symmetric spectra are necessary for this construction, and imposing more symmetry (for example orthogonal group actions) is not possible.  Multiplicativity of the unit morphism $\bbS_k \to \mtm$ required that the image of $M$ under each $e_m$ lie in a $\Sigma_m$-fixed subspace of $\mRk$.  These do not exist if, for instance, we allow the full orthogonal group to act.
\end{rmk}

As we have just seen, the symmetric sequence $\bbS_k$ is a useful model of the sphere spectrum for mapping to Thom spectra.  However, it is inconvenient for receiving collapse maps.  Luckily, our construction above supplies us with an abundance of multiplicative models for the sphere spectrum which are ideal targets for Thom collapses.  These will be important in the next section for comparing our construction with classical Atiyah duality.

To this end, observe that the closed $\epsilon$-neighborhood of $0 \in \bbR^k$ is a tubular neighborhood for the embedding $e \colon pt \into \bbR^k$ it defines, and the above construction gives a commutative $\bbS_k$-algebra corresponding to this data.  We write $\bbS_{k,\epsilon}$ for this algebra.

\begin{cor}
The unit map $\bbS_k \to \bbS_{k,\epsilon}$ is a level weak equivalence of commutative $\bbS_k$-algebras.
\end{cor}

An immediate first use of this construction is in defining augmentations on our ring spectra $\mtm$.  To accomplish this, we make the following more general observation.

\begin{lemma}
Let $N \subseteq M$ be a submanifold, and suppose $\epsilon < \min\{L_e, L_{e|_N}\}$.  Then there is a collapse map $\mtm \to \ntn$ compatible with the $\bbS_k$-algebra structures.
\begin{proof}
We use subscript $M$ and $N$ to distinguish between the symmetric sequences $D$ and $S$ for the two manifolds.  Note that $D_N (m)$ is a $\Sigma_m$-invariant subspace of $D_M (m)$ for every $m$, and similarly for $S_N (m)$ and $S_M (m)$.  We also have the inclusion
\begin{equation*}
S_M (m) \cap D_N (m) \subseteq S_N (m).
\end{equation*}
There is thus an induced morphism of symmetric sequences $\mtm \to \ntn$.  Compatibility with the unit morphisms from $\bbS_k$ may be read off from the definition of this map.  That the collapse is multiplicative is more tedious but no more difficult to verify.
\end{proof}
\end{lemma}

Now suppose $M$ is given a basepoint $x$.  In this case we require $e(x) = 0$.  The collapse map above now takes the form of a $\bbS_k$-algebra morphism $\mtm \to \bbS_{k,\epsilon}$.  This is not an augmentation on the point-set level, but becomes one once we pass to $\infty$-categories, where we are allowed to invert the weak equivalence of Corollary 3.1.

\subsection{Structured Atiyah Duality}
Recall that the Atiyah duality equivalence \cite{a61} is given as the adjoint of a map $\tilde \alpha \colon M_+ \wedge M^{\eta_e} \to S^k$ defined by
\begin{equation*}
\tilde \alpha (x, v) = v - e(x) \in \bbR^k / (\bbR^k \setminus B_\epsilon (0)) \cong \overline{B}_\epsilon (0) / \p \overline{B}_\epsilon (0) \cong S^k.
\end{equation*}
Classically, one shows that $\alpha$ is a stable duality in that the induced map
\begin{equation*}
\Sigma^{\infty-k} M^{\eta_e} \to F(M_+, \bbS)
\end{equation*}
is a weak equivalence of prespectra.  Our next result is a highly-structured version of this equivalence.  Note that since $M$ is homeomorphic to a CW complex, any level equivalence of symmetric sequences $A \to B$ induces a level equivalence on the cotensors $F(M_+, A) \to F(M_+, B)$.  In particular, $F(M_+, \bbS_k) \to F(M_+, \bbS_{k,\epsilon})$ is an equivalence of commutative $\bbS_k$-algebras, so as far as the homotopy theory of $E_\infty$-rings is concerned they are indistinguishable.

We now define the structured Atiyah duality morphism on the level of symmetric sequences.  Let $\tilde \alpha(m) \colon M_+ \wedge \mtm (m) \to \bbS_{k,\epsilon} (m)$ be given by the formula
\begin{equation*}
\tilde \alpha(m)(x, v) = v - e_m(x) \in \mRk / (\mRk \setminus B_\epsilon (0)) \cong \overline{B}_\epsilon (0) / \p \overline{B}_\epsilon (0) = \bbS_{k,\epsilon} (m).
\end{equation*}
By adjunction, this gives a map
\begin{equation*}
\alpha(m) \colon \mtm(m) \to F(M_+, \bbS_{k,\epsilon} (m)).
\end{equation*}
Equivariance is a consequence of all subsets in sight being invariant for the relevant symmetric group actions.

\begin{lemma}
The morphism $\alpha \colon \mtm \to F(M_+, \bbS_{k,\epsilon})$ defined above is compatible with the multiplications on the source and target.
\begin{proof}
This is a consequence of the formula $e_{m+n} (x) = e_m (x) + e_n (x) \in \mnRk$ for all $m,n$.
\end{proof}
\end{lemma}

The map $\alpha$ does not strictly commute with the unit maps for the algebra structures.  However, one can construct a homotopy-coherent triangle of multiplicative maps of the following form.
\begin{equation*}
\begin{tikzcd}
& \bbS_k \arrow{dl}[swap]{\eta} \arrow{dr}{c} \\
\mtm \arrow{rr}[swap]{\alpha} & & F(M_+, \bbS_{k,\epsilon})
\end{tikzcd}
\end{equation*}
More precisely, we define a homotopy $h_t$ from $\alpha \circ \eta$ to $c$ levelwise by
\begin{equation*}
h_t (v) (x) = [v - t \cdot e_m(x)] \in \bbS_{k,\epsilon} (m).
\end{equation*}
Equivariance is a consequence of the fact that $e_m(M)$ lies in a $\Sigma_m$-invariant linear subspace of $\mRk$ for all $m$.  As it stands, we do not have a morphism of symmetric spectra.  We resolve this in two steps, beginning with the following.

\begin{lemma}
Consider the symmetric sequence $F(M_+, \bbS_{k,\epsilon})$ as an $\bbS_k$-module via $\alpha \circ \eta$ rather than $c$.  Then $\alpha$ is a $\pi_*$-isomorphism of $k$-symmetric spectra.
\begin{proof}
By \cite{a61}, the map $\alpha(m)$ is a stable equivalence of spaces for each $m$.  Since $\mtm(m)$ and $F(M_+, \bbS_{k,\epsilon} (m))$ are both $(km-d)$-connected, $\alpha(m)$ is an equivalence on homotopy groups in dimensions $< 2(km-d)$ by Freudenthal's theorem.  By cofinality, $\alpha$ is a $\pi_*$-isomorphism.
\end{proof}
\end{lemma}

In order to complete the comparison, we employ coherence machinery due to Lurie in the context of algebra objects in monoidal $\infty$-categories.  In the previous Lemma, we considered $F(M_+, \bbS_{k,\epsilon})$ as a quasi-unital $E_\infty$-ring.  That is, the unit property of $c$ does not hold on the nose, merely up to a specified homotopy.  The primary result we need about these objects is the following.

\begin{thm}[{\cite[5.2.3.12]{ha}}]
The forgetful functor from the $\infty$-category of $E_\infty$-rings to the $\infty$-category of quasi-unital $E_\infty$-rings is an equivalence.
\end{thm}

The content of this result is that if a unit exists for a commutative monoid in an $\infty$-category, it is unique up to a contractible space of choices.  This has the following consequences.

\begin{cor}
The map $\alpha$ induces an equivalence in the $\infty$-category of $E_\infty$-rings $\mtm \to F(M_+, \bbS)$.  In particular, $\mtm$ has the correct homotopy type independent of the pair $(e, \epsilon)$.
\end{cor}

\begin{cor}
There is a zig-zag of weak equivalences of commutative $k$-symmetric ring spectra between $\mtm$ and $F(M_+, \bbS_k)$.
\end{cor}

The latter assertion is a result of the fact that we have a presentation of the $\infty$-category of $E_\infty$-rings as the model category of commutative $k$-symmetric ring spectra.

Similarly, we have seen that the choice of a point $x \in M$ determines ``augmentations'' $\mtm \to \bbS_{k,\epsilon}$ and $F(M_+, \bbS_{k,\epsilon}) \to \bbS_{k,\epsilon}$, and $\alpha$ strictly commutes with these morphisms.  Using the equivalence of Corollary 3.3, we may use these to view the two $E_\infty$-rings as augmented $E_\infty$-rings, and we obtain the following result, after making similar but easier coherence arguments as were necessary for the units.

\begin{prop}
The map $\alpha$ induces an equivalence in the $\infty$-category of augmented $E_\infty$ rings.
\end{prop}

\subsection{Relation to Previous Work}

In \cite{c04}, two related but larger models are given for $\mtm$.  In the largest, one tracks all $\epsilon$-tubular neighborhoods of all embeddings of $M$ into all euclidean spaces.  This model has great functoriality properties, but unfortunately does not have a strict unit given by the Thom collapse.  To correct this, the second model begins with a fixed embedding and $\epsilon$-neighborhood, and then tracks all linear isometric embeddings into higher-dimensional euclidean spaces.  Our model is actually a subspectrum of this latter one, obtained by noting that the diagonal embeddings are fixed by the symmetric group actions.

\begin{rmk}
Our construction gives a small, highly structured model for the Atiyah dual of a {\em single} manifold $M$, but does not appear to have good functoriality properties.  Given a map of manifolds $f \colon N \to M$, we may construct a structured umkehr or ``wrong-way'' map $f^! \colon \mtm \to \ntn$ if we allow the model of $\ntn$ to depend on $f$.  In the case that $f$ is an embedding of a submanifold, we have already seen a version of this construction in Lemma 3.2.  In general, it is constructed in essentially the same manner as in \cite[\S 4]{gs08}, with the relatively small addition of coherence machinery tracking higher diagonal embeddings.  We omit a detailed description since it is irrelevant to the remainder of our work.
\end{rmk}

\section{The Extended Cohen-Jones Construction}
One can understand the ring structure on $\mtm$ above in the following manner.  Given a manifold $M$, one can perform a Thom collapse on a tubular neighborhood of the diagonal embedding to induce a ``twisted'' product map $M \times M \to M^{TM}$.  To obtain an honest product, one ``untwists'' the whole construction by the stable normal bundle $-TM$.  This is expressed via the formula
\begin{equation*}
\mtm \wedge \mtm \to M^{-2TM + TM} \simeq \mtm.
\end{equation*}

Now consider a homotopy pullback diagram
\begin{equation*}
\begin{tikzcd}
Y \arrow{r} \arrow{d} & X \arrow{d}{p} \\
N \arrow{r}[swap]{e} & M
\end{tikzcd}
\end{equation*}
where the bottom arrow is an embedding of closed finite-dimensional manifolds.  Assume for now that $p$ is actually a fiber bundle.  A fundamental observation of Cohen and Klein in \cite{ck09} (extending the construction in \cite{cj02}) is that the upper arrow is now also an embedding of finite codimension, and moreover has a regular neighborhood isomorphic to the pulled-back normal bundle $p^* \eta_e$.  This good point-set behavior allows the definition of a pulled-back collapse map $X \to Y^{p^* \eta_e}$.  Moreover, one can twist this construction by a virtual bundle just like before.

With this in mind, we consider a fiber bundle $p \colon X \to M$, and perform the construction to the following pullback diagram.
\begin{equation*}
\begin{tikzcd}
X \times_M X \arrow{r}{\tilde \Delta} \arrow{d}[swap]{q} & X \times X \arrow{d}{p \times p} \\
M \arrow{r}[swap]{\Delta} & M \times M
\end{tikzcd}
\end{equation*}
One sees that the pullback of a tubular neighborhood of $\Delta(M)$ along $p \times p$ is a neighborhood of $\tilde \Delta(X \times_M X)$ homeomorphic to the pullback of $TM$ along $q$.  The collapse map from above gives a map $X \times X \to (X \times_M X)^{q^*TM}$.  In the case that $X$ is a fiberwise monoid, one may then compose with a map $(X \times_M X)^{q^*TM} \to X^{\, p^*TM}$ induced by multiplication.  This was done in the case of $X = LM$ in \cite{cj02}, and for general fiberwise monoids in \cite{gs08}.  After twisting by the virtual bundle $p^*(-TM)$ one obtains a spectrum $X^{-TM}$ which is a ring spectrum in the classical sense, i.e. a monoid in the stable category.

In the general case we consider, $X$ may be an arbitrary space over $M$, and so regular neighborhoods may not be well-behaved.  Instead, we focus on the disk bundles (or, ultimately, spherical fibrations) corresponding to the pulled-back stable normal bundles.  This turns out to be much better-behaved.  In particular, it preserves all weak equivalences between parametrized spaces, as we will now see.

\section{Generalized String Topology}
We are now prepared to give our main construction.  As in Section 3, fix an embedding $e \colon M \into \bbR^k$ and an $\epsilon \in (0, L_e)$.  We will suppress the pair $(e, \epsilon)$ from the notation, but every result stated in this section should be taken relative to this choice. We will give a universal characterization of this functor in Section 8 using Morita theory in $\infty$-categories.  Also recall the definition of the symmetric sequences $D$ and $S$ from Section 3, our models of disk and sphere bundles in the stable normal bundle of $M$, respectively.

The generalized string topology functor based on $(e, \epsilon)$ is a composite of four functors.  We will show that each step preserves (point-set) colimits and tensors with (unbased) spaces.

\begin{const}
Assign to a space $X \in \sT_{/M}$ the constant symmetric sequence $\bar X$ in $\sT_{/M}$, i.e. $\bar X(m) = X$ with the trivial $\Sigma_m$-action.

\begin{lemma}
The functor $X \mapsto \bar X$ is a left adjoint and preserves tensors with spaces.
\begin{proof}
This is just the forgetful functor associated to the continuous functor of topological categories $\Sigma \to 0$, and thus has all these properties from elementary enriched category theory.
\end{proof}
\end{lemma}

Next, pull back the symmetric sequence $\bar X$ over the levelwise Hurewicz cofibration of symmetric sequences $S \to D$, i.e. apply to $\bar X$ the functor $i_{q \to p}^* \colon \Sigma \sT_{/M} \to Ar(\Sigma \sT)_{/ S \to D}$ given by the formula
\begin{equation*}
\{E(m)\} \mapsto \{E(m) \times_M S(m)\} \to \{E(m) \times_M D(m)\}.
\end{equation*}
We introduce the shorthand $S_X \to D_X$ for the image of $\bar X$ under $i_{q \to p}^*$.  Note that by \cite{kie87}, this resulting map of symmetric sequences is also a levelwise Hurewicz cofibration.

\begin{lemma}
The pullback functor $i_{q \to p}^* \colon \Sigma \sT_{/M} \to Ar(\Sigma \sT)_{/S\to D}$ preserves colimits and tensors with spaces.
\begin{proof}
Suppose $f \colon A \to B$ is a continuous map of compactly generated spaces.  Then the pullback along $f$ induces a functor $f^* \colon \sT_{/B} \to \sT_{/A}$.  In \cite[2.1]{ms06}, the authors construct a right adjoint $f_*$ to this functor.  Essentially, it takes a space $X \to A$, viewed as a space over $B$ by composing with $f$, and assigns to it a space of relative sections $A \to X$.  There are point-set topological difficulties in general, but if $f$ is an open map (for example a fiber bundle whose fibers are CW complexes) then there are no issues.  Since $S(m)$ and $D(m)$ are fiber bundles over $M$ with smooth manifold fibers, the relevant right adjoint functors exist.

Now let $X \colon I \to \Sigma \sT_{/M}$ be a small diagram.  Recall that colimits in diagram categories are computed pointwise, i.e. we have a natural isomorphism
\begin{equation*}
(\colim_I X)(m) \cong \colim_I (X(m)).
\end{equation*}
A corresponding formula holds in $Ar(\Sigma \sT)_{/S \to D}$, and so to show preservation of colimits we need only check that the resulting cone on $i_{q \to p}^* \circ X$ is a colimit cone at each object in $[1] \times \Sigma$.  At $(0,m) \in [1] \times \Sigma$, the diagram takes the form
\begin{equation*}
i \mapsto X(i) \times_M S(m),
\end{equation*}
which is equivalent to applying $- \times_M S(m)$ to the diagram $X$ restricted to the object $m$.  But we know this functor preserves colimits by our earlier construction, so
\begin{equation*}
(\colim_I (i_{q \to p}^* \circ X))(0,m) \cong \colim_i (i_{q \to p}^* \circ X)(0,m).
\end{equation*}
A similar argument applies to $(1,m) \in [1] \times \Sigma$, and so $i_{q \to p}^*$ preserves colimits.

The proof that $i_{q \to p}^*$ preserves tensors is analogous but easier, the main point being that tensors in diagrams of spaces are computed pointwise.
\end{proof}
\end{lemma}

Next, to the morphism of symmetric sequences $S_X \to D_X$ assign the symmetric sequence $\{D_X (m)/S_X(m)\}$ over $\mtm$,  i.e. apply the functor from $Ar(\Sigma \sT)_{/S \to D}$ to $(\Sigma \sT_*)_{/\mtm}$ that assigns to a map $g \colon E \to F$ of symmetric sequences over the morphism $S \to D$ levelwise quotient symmetric sequence $E/F$ over $D/S = \mtm$.  Write $\xtm$ for the result.

\begin{lemma}
The levelwise quotient functor $Ar(\Sigma \sT)_{/S \to D} \to (\Sigma \sT_*)_{/\mtm}$ preserves colimits and tensors with spaces.
\begin{proof}
This functor is itself a colimit functor, and colimits commute.  Preservation of tensors follows from the identification
\begin{equation*}
(X \times Y)/(A \times Y) \cong (X/A) \wedge Y_+
\end{equation*}
for all pairs of spaces $(X,A)$ and all spaces $Y$.
\end{proof}
\end{lemma}

Finally, forget the morphism to $\mtm$ to obtain a symmetric sequence $\xtm$.

\begin{lemma}
The forgetful functor $(\Sigma \sT_*)_{/\mtm} \to \Sigma \sT_*$ is a left adjoint and preserves tensors with spaces.
\begin{proof}
The right adjoint functor $\Sigma \sT_* \to (\Sigma \sT_*)_{/\mtm}$ is given by taking the categorical product with $\mtm$ in $\Sigma \sT_*$.  Preservation of tensors follows from the fact that they are computed in $\Sigma \sT_*$.
\end{proof}
\end{lemma}
\end{const}

Combining the above lemmas, we have the following.

\begin{prop}
The generalized string topology functor $(X \to M) \mapsto \xtm$ is a continuous left adjoint and preserves tensors with spaces, when viewed as a functor into $\Sigma \sT_*$.
\end{prop}

The above construction is not so useful as it stands.  We would like $(-)^{-TM}$ to land in a category of modules for the Atiyah dual $\mtm$.  The following shows that this is, in fact, the case.  Recall that a functor $F \colon \sC \to \sD$ between symmetric monoidal categories is lax symmetric monoidal if there is a morphism $\mathbbm{1}_\sD \to F(\mathbbm{1}_\sC)$ and a bi-natural transformation $F (-) \otimes_\sD F (-) \Rightarrow F(- \otimes_\sC -)$ satisfying natural compatibility properties with the associativity and commutativity transformations for the symmetric monoidal structures.  See \cite[\S 20]{mmss} for a more thorough discussion.

\begin{lemma}
The generalized string topology functor is lax symmetric monoidal as a functor to symmetric sequences.
\begin{proof}
Suppose $X, Y \in \sT_{/M}$.  We must produce natural maps
\begin{equation*}
\mu \colon \xtm (m) \wedge \ytm (n) \to \xytm (m+n)
\end{equation*}
for all natural numbers $m,n$ which are appropriately equivariant, associative and commutative.  To this end, we observe the Hurewicz cofibrations
\begin{equation*}
D_{X \times_M Y} (m+n) \subset D_X (m) \times D_Y (n) \subset X \times Y \times \mnRk
\end{equation*}
and
\begin{equation*}
D_{X \times_M Y} (m+n) \cap ((S_X (m) \times D_Y (n)) \cup (D_X (m) \times S_Y (n))) \subset S_{X \times_M Y} (m+n).
\end{equation*}
It follows that we have a collapse map $\mu$ as claimed.  All basic properties of this map may be verified by easy but tedious diagram chases involving the action of $\Sigma_m$ on $\mRk$ and coherence properties of the symmetric monoidal products $\wedge$ and $\times_M$.  The most subtle is perhaps compatibility with the twist isomorphisms, so we include this case as an example and leave the rest to the reader.  We wish to show that the following diagram commutes.
\begin{equation*}
\begin{tikzcd}
\xtm \wedge \ytm \arrow{r} \arrow{d}[swap]{\tau} & (X \times_M Y)^{-TM} \arrow{d}{\tau^{-TM}} \\
\ytm \wedge \xtm \arrow{r} & (Y \times_M X)^{-TM}
\end{tikzcd}
\end{equation*}
Restricting to the wedge factor corresponding to $\xtm(m) \wedge \ytm(n)$ in the upper-left corner, we recall that the twist isomorphism in symmetric sequences includes an action by the block transposition $\chi_{m,n}$, which then acts on $(Y \times X)^{-TM}$ by permuting coordinates in $\mnRk$.  Hence, even though the left-hand arrow acts on $-TM$ while the right-hand arrow does not, the symmetric sequence structure on everything in sight forces the bottom arrow to ``undo'' this action.
\end{proof}
\end{lemma}

Since $\sT_{/M}$ is considered with the Cartesian monoidal structure, every space over $M$ is naturally a module over the final object $M$, considered as a commutative algebra via the unit isomorphism $M \times_M M \cong M$.  The lax symmetric monoidal structure of $(-)^{-TM}$ thus has the following consequence.

\begin{thm}
The generalized string topology functor lifts to a lax symmetric monoidal functor $(-)^{-TM} \colon \sT_{/M} \to \Sigma \sS_{\mtm}$ which preserves weak equivalences.
\begin{proof}
First, we show that for any space $X \in \sT_{/M}$, the symmetric sequence $\xtm$ is naturally a module for $\mtm$.  We have a natural isomorphism $X \times_M M \cong X$ in spaces over $M$.  Thus there are natural morphisms of symmetric sequences
\begin{equation*}
\xtm \wedge \mtm \to (X \times_M M)^{-TM} \cong \xtm.
\end{equation*}
The commutativity of $\mtm$ then gives $\xtm$ the structure of a left module as well as a right module.  Using the associativity and unit transformations for the symmetric monoidal structures on $\Sigma \sT$ and $\sT_{/M}$ along with the natural isomorphism
\begin{equation*}
X \times_M M \times_M Y \cong X \times_M Y,
\end{equation*}
we see that the transformation
\begin{equation*}
\xtm \wedge \ytm \to (X \times_M Y)^{-TM}
\end{equation*}
factors through the quotient
\begin{equation*}
\xtm \wedge \ytm \to \xtm \wedge_{\mtm} \ytm.
\end{equation*}
From here, and a few diagram chases, it follows that $(-)^{-TM}$ is lax symmetric monoidal as a functor to $\Sigma \sS_{\mtm}$.

We will now show that generalized string topology preserves weak equivalences.  Suppose $f \colon X \to Y$ is a weak equivalence of spaces over $M$.  By the five lemma applied to the long exact sequence of a fibration, the induced maps $D_X (m) \to D_Y (m)$ and $S_X (m) \to S_Y (m)$ are weak equivalences for all $m$.  Recall that for each $m$ and any $Z \to M$, the maps $S_Z (m) \into D_Z (m)$ are Hurewicz cofibrations, and so 
\begin{equation*}
S_Z (m) \into D_Z (m) \to Z^{-TM} (m)
\end{equation*}
is a cofiber sequence.  We may now apply five lemma to the long exact sequences for homology to conclude that $\xtm (m) \to \ytm (m)$ is a homology equivalence for each $m$.  If $m > 1$, these spaces are simply connected, and therefore the map $\xtm \to \ytm$ is a $\pi_*$-equivalence, in particular a stable equivalence.  Since weak equivalences in $\mtm$-modules are created in $\bbS_k$-modules, $(-)^{-TM}$ preserves weak equivalences.
\end{proof}
\end{thm}

\begin{rmk}
There is a natural version of this functor for ex-spaces over $M$.  All the previous arguments apply equally well to this case, except that $(-)^{-TM}$ will in general only preserve weak equivalences between ex-spaces for which the section is a Hurewicz cofibration.  This is the first step towards constructing a generalized string topology functor for spectra parametrized by $M$.  We omit details, as this extension will be handled differently in Section 6.
\end{rmk}

We conclude this section with an immediate corollary.  Fix a connected, based CW complex $X$.  Then using the symmetric monoidal structure $\times_X$, there is a natural notion of an operad $O$ in spaces over $X$, and of an $O$-algebra $A$ over $X$.  We call such an $O$ a fiberwise operad, and such an $A$ a fiberwise $O$-algebra.  One immediate source of fiberwise operads is the usual category of operads in topological spaces.  Namely, given an operad $O$, we obtain a fiberwise operad $O \times X$ by taking the levelwise product with $X$.  It turns out that the action of the fiberwise operad $O \times X$ on a space $A$ over $X$ is equivalent to the action of $O$ on $A$ via the tensoring of $\sT_{/X}$ over $\sT$.  In either case, we abusively call $A$ a fiberwise $O$-algebra.  This is among the first cases considered in the string topology literature.  In particular, the space of (unbased) maps $S^n \to X$ is a fiberwise $E_n$-algebra (essentially because each fiber is isomorphic to the based $n$-fold loop space of $X$.)

\begin{cor}
Suppose $O$ is an operad in spaces and $A \to M$ is a fiberwise $O$-algebra.  Then $A^{-TM}$ is an $O$-algebra in $\mtm$-modules.  In particular, if $A$ is a fiberwise $E_n$-algebra, then $A^{-TM}$ is an $E_n$-ring spectrum.
\begin{proof}
Generalized string topology is lax symmetric monoidal and preserves tensors with spaces.
\end{proof}
\end{cor}

\section{Homotopy Colimits}

Our characterization of generalized string topology in the next section requires an analysis of the behavior of $(-)^{-TM}$ with respect to homotopy colimits.  For this, we employ the theory developed in \cite{shu06}.  In that work, Shulman gives a comparison between local and global definitions of homotopy colimits.  Here ``local'' refers to using a simplicial enrichment to explicitly construct homotopy coherent cones, and ``global'' means using model category-theoretic framework to deform the colimit functor to one that preserves weak equivalences.  The latter approach is very well-suited to situations where one has good control over model-theoretic cofibrations.  The former is not so reliant on this kind of structure, and instead the most delicate part of the construction is comparing the homotopy theory coming from the enrichment with that coming from the model structure.  We first recall the relevant definitions, and then apply the theory to our case.

\begin{defn}
Let $I$ be a small category, $\sC$ a simplicial (or topological) model category, $X \colon I \to \sC$ a functor.  The uncorrected homotopy colimit of $X$ is the coend
\begin{equation*}
\uhocolim_I X = \int^{i \in I} N(I_{i/}) \otimes X(i).
\end{equation*}
The corrected homotopy colimit of $X$ is obtained by applying $\hocolim$ to an objectwise cofibrant replacement of $X$, i.e.
\begin{equation*}
\hocolim_I X = \uhocolim_I (Q \circ X).
\end{equation*}
\end{defn}

It is shown in \cite[\S 9]{shu06} that this definition of homotopy colimit agrees with the left-derived functor of $\colim \colon \sC^I \to \sC$, and so coincides with the usual notion.

\begin{rmk}
Unraveling the definitions, it is the case that the coend in the definition above is precisely the geometric realization of a simplicial prespectrum given by the bar construction and written $B_\bullet (*, I, X)$.  It has the additional property that all of the degeneracy morphisms are Hurewicz cofibrations.  This means that homotopy colimits in $\sP_k$ behave well with respect to homotopy equivalences, as we will see below.
\end{rmk}

\begin{lemma}
Generalized string topology preserves uncorrected homotopy colimits.
\begin{proof}
This follows from preservation of tensors and ordinary colimits.
\end{proof}
\end{lemma}

Unfortunately, $(-)^{-TM}$ appears not to send CW complexes over $M$ to cofibrant $\mtm$-modules, so pointwise cofibrant replacement in spaces over $M$ does not obviously compute the correct homotopy colimit.  With the help of a little classical homotopy theory, though, we are saved.  First, we forget the $\mtm$-module structures and $\Sigma_m$-actions on the spectra $\xtm$ to view $(-)^{-TM}$ as a functor to the category of $k$-prespectra, which we denote by $\sP_k$.

\begin{lemma}
The functor $(-)^{-TM} \colon \sT_{/M} \to \sP_k$ preserves colimits and tensors with spaces, and so preserves uncorrected homotopy colimits.
\begin{proof}
The forgetful functor $\Sigma \sS_k \to \sP_k$ is associated to the topological functor $\sN_{\bbS_k} \to \Sigma_{\bbS_k}$, and so preserves colimits and tensors.
\end{proof}
\end{lemma}

We now show that this composite functor sends cofibrant spaces over $M$ to $k$-prespectra homotopy equivalent to cofibrant $k$-prespectra.  This will imply the preservation of corrected homotopy colimits.  First we need the following definition.

\begin{defn}
We say $E \in \sP_k$ is $\Sigma^k$-cofibrant if the structure maps $E(m) \wedge S^k \to E(m+1)$ are Hurewicz cofibrations for all $m$.  Let $K \colon \sP_k \to \sP_k$ be the functor defined levelwise by
\begin{equation*}
KE(m) = \operatorname{tel}_{n \leq m} E(n) \wedge S^{m-n},
\end{equation*}
where the right-hand side is the mapping telescope built from the iterated structure maps.  The structure maps on $KE$ are the natural ones.
\end{defn}

One can quickly check that $KE$ is $\Sigma^k$-cofibrant for all $E$.  There is also a natural level equivalence $K \Rightarrow \id$ given by ``collapsing the telescope'' though the structure maps for $E$.

\begin{lemma}
The functor $K$ preserves colimits and tensors with spaces.
\begin{proof}
Colimits and tensors with spaces commute, and $K$ is constructed as a combination of those operations.
\end{proof}
\end{lemma}

We now recall the following variant of a result from \cite{mmss}, whose proof works equally well in our case.

\begin{prop}[{\cite[11.4]{mmss}}]
If $E \in \sP_k$ is $\Sigma^k$-cofibrant and every space $E(m)$ has the homotopy type of a CW complex, then $E$ has the homotopy type of a cofibrant $k$-prespectrum.
\end{prop}

\begin{cor}
If $X \in \sT_{/M}$ is cofibrant, then $K\xtm \in \sP_k$ has the homotopy type of a cofibrant $k$-prespectrum.
\begin{proof}
Since $X$ is cofibrant, it has the homotopy type of a CW complex.  Each space $D_X (m)$ and $S_X (m)$ thus has the homotopy type of a CW complex, and therefore so does $\xtm(m)$.  From the definition of $K$ we see that each $K\xtm(m)$ has the homotopy type of a CW complex, and so we may apply the previous Proposition.
\end{proof}
\end{cor}

The following lemma shows that in the situation above, the uncorrected homotopy colimit of $KX$ is equivalent to the corrected homotopy colimit.

\begin{lemma}
Suppose $\alpha \colon X \Rightarrow X'$ is a natural transformation of functors $I \to \sP_k$ which is an objectwise homotopy equivalence.  Then the induced map
\begin{equation*}
\uhocolim_I X \to \uhocolim_I X'
\end{equation*}
is a homotopy equivalence.  In particular, if $X$ is as in Corollary 6.8, then
\begin{equation*}
\uhocolim_I K\xtm \simeq \hocolim_I K\xtm.
\end{equation*}
\begin{proof}
Recall from Remark 6.2 that the uncorrected homotopy colimits are geometric realizations of simplicial $k$-prespectra all of whose degeneracies are Hurewicz cofibrations.  Moreover, the induced simplicial map
\begin{equation*}
B_\bullet(*,I,X) \to B_\bullet(*,I,X')
\end{equation*}
is a levelwise homotopy equivalence.  The induction argument given for \cite[X.2.4]{ekmm} may be easily adapted to our case to prove that the induced map of geometric realizations is a homotopy equivalence as claimed.  For the last statement, note that the natural transformation $QKX^{-TM} \Rightarrow KX^{-TM}$ is an objectwise homotopy equivalence by 2-of-3 for homotopy equivalences.
\end{proof}
\end{lemma}

Piecing it all together, we obtain our desired result.

\begin{prop}
The generalized string topology functor $(-)^{-TM} \colon \sT_{/M} \to \sP_k$ preserves homotopy colimits up to weak equivalence.
\begin{proof}
First, we restrict to CW complexes over $M$ so as to compute the correct homotopy colimits in the domain.  Let $X \colon I \to \sT_{/M}$ be a diagram of CW complexes over $M$.  Then
\begin{align*}
(\hocolim_I X)^{-TM} & \simeq (\uhocolim_I X)^{-TM} \\
& \simeq K (\uhocolim_I X)^{-TM} \\
& \cong \uhocolim_I K\xtm \\
& \simeq \hocolim_I K\xtm \\
& \simeq \hocolim_I \xtm.
\end{align*}
Only the last equivalence has not already been argued, but it follows from the fact that the corrected homotopy colimit functor respects objectwise weak equivalences of diagrams.
\end{proof}
\end{prop}

This implies the following, which is the last required component for our $\infty$-categorical analysis.

\begin{cor}
The induced functor $(-)^{-TM} \colon \cT_{/M} \to \Mod_{\mtm}$ preserves ($\infty$-categorical) colimits.
\begin{proof}
Since $\mtm \wedge -$ preserves colimits in the $\infty$-category of spectra, colimits in $\Mod_{\mtm}$ are created in the $\infty$-category $\cS$ of spectra.  Since $\sP_k$ models $\cS$,  \cite[4.2.4.1]{ha} implies that the functor underlying generalized string topology preserves colimits.
\end{proof}
\end{cor}

\subsection{Extension to Parametrized Spectra}

The previous corollary allows us to use the theory of $\infty$-categories to extend generalized string topology to parametrized spectra while avoiding the perils of point-set topology.  More precisely, applying Proposition 2.11 we see that $(-)^{-TM}$ extends to the $\infty$-category of spectrum objects in spaces over $M$.  One of the great virtues of the $\infty$-categorical context is the following theorem.

\begin{thm}[{\cite[3.6]{abg10}}]
The stable $\infty$-category $\Sp(\cT_{/M})$ is equivalent to the $\infty$-category underlying the stable model category $\sS_{M}$ of orthogonal spectra parametrized by $M$, as constructed in \cite[\S12.3]{ms06}.  This equivalence respects the symmetric monoidal structure $\wedge_M$, and the functor $\Sigma^\infty_+ \colon \cT_{/M} \to \Sp(\cT_{/M})$ corresponds to the fiberwise suspension spectrum functor.
\end{thm}

From now on, we write $\cS_{/M}$ for the stable $\infty$-category of spectra parametrized by $M$, and $\Sigma^\infty_M$ for the fiberwise suspension spectrum functor.  We may now immediately deduce the following.

\begin{cor}
The generalized string topology functor $(-)^{-TM}$ extends uniquely through $\Sigma^\infty_M$ to a colimit-preserving functor
\begin{equation*}
(-)^{-TM} \colon \cS_{/M} \to \Mod_{\mtm}
\end{equation*} 
\end{cor}

This induced functor is also lax symmetric monoidal, though we will defer the proof to a sequel in which we explore the construction of string topology operations on twisted generalized (co)homology of spaces over $M$.

\section{Derived Koszul Duality}

We will here review the versions of stable and unstable derived Koszul duality, alias bar-cobar duality, that we require.  First, consider the following unstable or ``nonabelian'' situation.  Let $X$ be a connected based CW-complex.  Then $X \simeq B\Omega X$, and there is a close relationship between spaces with an action of $\Omega X$ and spaces with a map to $X$.  To fit this with the usual notion of bar-cobar duality, we remark that there is an isomorphism of categories between spaces over $X$ and comodules for the coalgebra $(X, \Delta)$.  Therefore we are comparing the categories of modules for a monoid $\Omega X$ and comodules for its bar construction $X$.  Let $PX$ be the based path space of $X$.  Then, assuming appropriate models are chosen, the contractible space $PX$ comes equipped with a fibration $PX \to X$ and a free right action by $\Omega X$, and so plays the role of the universal principal $\Omega X$-bundle $E \Omega X \to B \Omega X$.  For concreteness, one may take $\Omega X$ to be the (strictly associative) Moore loop space of $X$, and $P X$ the Moore path space of $X$.

\begin{prop}[{\cite[8.5]{shu08}}]
Let $X$ be as above.  Then there is a Quillen equivalence
\begin{equation*}
(-) \times_{\Omega X}  PX \colon \sT^{\, \Omega X} \rightleftarrows \sT_{/X} \colon F_{X} (PX, -).
\end{equation*}
Here we consider each category with the $m$-model structures created through the associated forgetful functors to $\sT$.
\end{prop}

The form of the right adjoint functor is dictated by point-set considerations.  When restricted to cofibrant-fibrant spaces over $X$, there is a natural weak equivalence from $F_X (PX, -)$ to $PX \times_X (-)$, the usual homotopy fiber functor.  Thus, we may equally take (the derived functor of) the latter as our inverse equivalence to the induced functor of $\infty$-categories
\begin{equation*}
(-) \times_{\Omega X} * \colon \cT^{\, \Omega X} \to \cT_{/X}.
\end{equation*}
This need not really concern us, since we mostly make use of the above functor, rather than its inverse.

We now discuss the stable version of derived Koszul duality.  Let $k$ be an $E_\infty$-ring, $A$ an augmented $k$-algebra.

\begin{defn}
The Koszul dual of $A$ is is the derived endomorphism $A_\infty$-ring
\begin{equation*}
A^! := \End_A (k).
\end{equation*}
\end{defn}

The $A_\infty$-ring $A^!$ is naturally an augmented $k$-algebra, and there is a canonical map $A \to (A^!)^!$.  If it is an equivalence, we say $A$ is dc-complete.  To make contact with bar-cobar duality, we note that there is a natural equivalence of $A_\infty$-rings
\begin{equation*}
\End_A (k) \simeq \Hom_k (k \wedge_A k, k).
\end{equation*}
That is, $A^!$ is the $k$-linear dual of the coalgebra obtained from $A$ via the bar construction.

\begin{ex}
Suppose $X$ is a finite connected based CW-complex, and consider the $A_\infty$-ring $A := \Sigma^\infty_+ \Omega X$.  Then $A^! \simeq \bbD X$.  Moreover, $A$ is dc-complete if and only if $X$ is simply-connected.  See \cite[\S3]{bm11} for a proof.
\end{ex}

Just as in the nonabelian situation above, there are adjoint functors
\begin{equation*}
(-) \wedge_A k \colon \Mod_A \rightleftarrows \Mod_{A^!} \colon \Hom_{A^!} (k, -).
\end{equation*}
However, they do not always form an equivalence of underlying $\infty$-categories.  Some smallness is required of the augmentation $A \to k$, and even then one often only has an equivalence to a full subcategory.  The specific situation of $A = \Sigma^\infty_+ \Omega X$ is treated in detail in \cite{bm11}, and the general theory is developed in \cite{dgi06}.

\section{Algebraic Generalized String Topology}

In Section 3, we saw that any particular geometric construction of the spectrum $\mtm$ was equivalent to the Spanier-Whitehead dual $F(M_+, \bbS)$ as $E_\infty$ rings.  In this section, we wish to extend this result to the whole of the generalized string topology construction.  The context for our work is a kind of Morita theory in stable $\infty$-categories.  Recall the classical Eilenberg-Watts theorem.

\begin{thm}[Eilenberg-Watts]
Suppose $R, S$ are rings and
\begin{equation*}
F \colon \Mod_R \to \Mod_S
\end{equation*}
a colimit-preserving functor.  Then there is an $R-S$ bimodule $M$ (unique up to isomorphism) and a natural isomorphism of functors
\begin{equation*}
F \cong (-) \otimes_R M.
\end{equation*}
\end{thm}
The proof proceeds by observing that any $R$-module is isomorphic to the colimit of a diagram of free $R$-modules of rank 1, and so $F$ is entirely determined by its value on $R$ itself.  Preservation of colimits implies that $F$ respects the additive enrichment on the domain and codomain categories, and so provides a morphism $R \to \End_S(F(R))$.  That is, $F(R)$ is naturally an $R-S$ bimodule.  This is the $M$ of the theorem.

Now suppose that $R,S$ are $A_\infty$-rings, and consider a colimit-preserving functor
\begin{equation*}
F \colon \Mod_R \to \Mod_S.
\end{equation*}
For any stable $\infty$-category $\cA$ and any objects $x, y \in \cA$, there is a naturally defined spectrum $\Hom_\cA (x,y)$ whose underlying infinite loop space is equivalent to the space of maps from $x$ to $y$, and colimit-preserving (or more generally exact) functors may be shown to repect the resulting spectral enrichments.  With this in hand, the argument above applies more or less verbatim to prove the $\infty$-categorical analogue of the Eilenberg-Watts theorem.

\begin{thm}
Suppose $R, S$ are $A_\infty$-rings and
\begin{equation*}
F \colon \Mod_R \to \Mod_S
\end{equation*}
a colimit-preserving functor.  Then there is an $R-S$ bimodule $M$ (unique up to equivalence) and a natural equivalence of functors
\begin{equation*}
F \simeq (-) \wedge_R M.
\end{equation*}
\end{thm}
We refer the reader to \cite[\S 8]{abghr} for precise details.

In previous sections we have seen that $\cT_{/M}$ is equivalent to the $\infty$-category $\cT^{\, \Omega M}$, and that $(-)^{-TM}$ factors through the stabilization of its domain category.  The following observation then places us squarely in the above situation.

\begin{lemma}
There is an equivalence of categories
\begin{equation*}
\Sp(\cT^{\, \Omega M}) \simeq \Mod_{\Sigma^\infty_+ \Omega M}.
\end{equation*}
\end{lemma}

So, by precomposing $(-)^{-TM}$ with nonabelian bar-cobar duality and then factoring through the stabilization, we obtain a colimit-preserving functor between $\infty$-categories of modules over $A_\infty$-rings.  The following result identifies this functor in the spirit of the Eilenberg-Watts theorem.

\begin{thm}
Fix a basepoint $x \in M$.  The following diagram of presentable $\infty$-categories and colimit-preserving functors commutes up to natural equivalence.
\begin{equation*}
\begin{tikzcd}
\cT^{\, \Omega M} \arrow{rr}{(-) \times_{\Omega M} *}[swap]{\simeq} \arrow{d}[swap]{\Sigma^\infty_+} & & \cT_{/M} \arrow{dd}{(-)^{-TM}} \\
\operatorname{Mod}_{\Sigma^\infty_+ \Omega M} \arrow{d}[swap]{(-) \wedge_{\Sigma^\infty_+ \Omega M} S^{-TM}} & &  \\
\operatorname{Mod}_{\bbD M} \arrow{rr}{\simeq}[swap]{\alpha^*} & & \operatorname{Mod}_{\mtm}
\end{tikzcd}
\end{equation*}
\begin{proof}
By Theorem 8.2, it suffices to identify the image of $\Omega M$ as a left $\Sigma^\infty_+ \Omega M$-module in $\Mod_{\bbD M}$.  Passing along the top-right composite results in the space $*^{-TM} \simeq S^{-d}$.  On the point-set level, we use the space $PM$ as our model for $*$.  There is a natural map $\Omega M \into \Aut_M (PM)$ since $PM$ is the universal principal $\Omega M$-fibration, which is then sent to a morphism $\Omega M \to \End_{\bbD M} (S^{-d})$.  We denote this bimodule $S^{-TM}$.
\end{proof}
\end{thm}

Though this theorem takes place in the realm of $\infty$-categories, we in fact have models for the relevant categories involved as coherent nerves of simplicially localized model categories.  Thus the following point-set level statement can be deduced quite easily.

\begin{cor}
For any two choices of $(e, \epsilon)$, the resulting generalized string topology functors are connected by a zigzag of natural weak equivalences, whenever this makes sense.
\end{cor}

\appendix
\section{$k$-Diagram Spectra}

The theory of model categories of diagram spectra was developed in detail in \cite{mmss}.  Briefly, given a symmetric monoidal topological category $\sD$, one obtains a symmetric monoidal category $\sD \sT_*$ of $\sD$-shaped diagrams in based spaces, and given a commutative monoid $S$ in $\sD \sT_*$ one defines a symmetric monoidal topological category of modules over $S$.  All these categories carry topological model structures, and a certain Bousfield localization of the category of $S$-modules is the category of $\sD$-spectra over $S$.  The paper \cite{mmss} studies various examples, all of which are Quillen equivalent and model the stable homotopy category.

None of the results in this appendix are surprising, and essentially all of the work is done in \cite{mmss}, from which we borrow definitions and notation heavily.  However, we were unable to find the precise statements we need in the literature.

\subsection{Definitions}

Let $\sD$ be a symmetric monoidal topological category, $k$ a nonnegative integer.  Write $i_k \colon \sD \to \sD$ for the $k^{\text{th}}$ power functor, i.e. $X \mapsto X^{\otimes k}$.  It follows from an easy formal computation that $i_k$ is strong symmetric monoidal.  By precomposition, $i_k$ induces a lax symmetric monoidal functor $\bbU_k \colon \sD \sT_* \to \sD \sT_*$.  Write $\bbP_k$ for its (strong symmetric monoidal) left adjoint, a proof of whose existence may be found in \cite[\S 23]{mmss}.

The relevant cases at hand are when $\sD$ is one of the categories $\mathscr{N}$, $\Sigma$ or $\mathscr{I}$ corresponding to prespectra, symmetric spectra or orthogonal spectra, respectively.  In these cases, $i_k$ is faithful, corresponding to ``block'' symmetries.  Let $\bbS$ be the standard sphere $\sD$-space.  We write $\bbS_k$ for $\bbU_k (\bbS)$.  Concretely, $\bbS_k (m) = S^k \wedge ... \wedge S^k \cong S^{mk}$ with $\sD(m,m)$ acting in blocks.  For each $k$, $\bbS_k$ is a monoid in $\sD \sT_*$, and is commutative if $\bbS$ is.  We call right $\bbS_k$-modules $k$-$\sD$-spectra, and denote the category of them by $\sD \sS_k$.  We obtain the following from combining Propositions 3.5-3.8 of \cite{mmss}.

\begin{prop}
The forgetful functor $\bbU_k \colon \sD \sT_* \to \sD \sT_*$ induces forgetful functor $\bbU_k \colon \sD \sS \to \sD \sS_k$, and $\bbP_k$ similarly extends to a left asjoint prolongation functor $\bbP_k \colon \sD \sS_k \to \sD \sS$.  If $\bbS$ is commutative, then $\bbU_k$ is lax symmetric monoidal and $\bbP_k$ is strong symmetric monoidal.
\end{prop}

We consider these categories of diagram spaces and diagram spectra as tensored over the category of based spaces by the formula
\begin{equation*}
(X \otimes E)(d) = X \wedge E(d),
\end{equation*}
and over the category of unbased spaces by addition of a disjoint basepoint.  These categories are also cotensored and topologically enriched, and these structures are suitably compatible, see \cite[\S 1]{mmss}.

The results of \cite[\S 2]{mmss} shows that the categories $\sD \sS_k$ are isomorphic to categories of diagram spaces for a new category $\sD_{\bbS_k}$ built from the monoid $\bbS_k$.  In the case that $\bbS_k$ is commutative, the category $\sD_{\bbS_k}$ is symmetric monoidal, and the isomorphism with $\sD \sS_k$ is strong symmetric monoidal.  This allows us to develop the model theory of categories of $k$-$\sD$-spectra simultaneously with the model theory of $\sD$-spaces.

\subsection{The Level Model Structure}

\begin{defn}
Let $f \colon X \to Y$ be a morphism in $\sD \sT_*$.  We say $f$ is a level equivalence if $f(d) \colon X(d) \to Y(d)$ is a weak equivalence of based spaces for every $d \in \sD$.  Similarly, we call $f$ a level fibration if each $f(d)$ is a Serre fibration, and say $f$ is an acyclic level fibration if it is a level fibration and a level equivalence.  A $q$-cofibration is a morphism in $\sD \sT_*$ which has the left lifting property with respect to all acyclic level fibrations.
\end{defn}

The following is proved in \cite[\S 6]{mmss} as Theorem 6.5.

\begin{thm}
The category $\sD \sT_*$ is a compactly generated proper topological model category with respect to the level equivalences, level fibrations, and q-cofibrations.
\end{thm}

Specializing to the category $\sD = \sD_{\bbS_k}$ we obtain a level model structure on the category $\sD \sS_k$.

\subsection{Homotopy Groups, Stable Equivalences and the Stable Model Structure}

\begin{defn}
Let $X \in \sD \sS_k$, where $\sD \in \{ \mathcal{N}, \Sigma, \mathcal{I}\}$.  The homotopy groups of $X$ are the abelian groups defined by
\begin{equation*}
\pi_n (X) = \colim_{m \to \infty} \, [S^{n + km}, X(m)].
\end{equation*}
A morphism $f \colon X \to Y$ in $\sD \sS_k$ is a $\pi_*$-isomorphism if $f$ induces isomorphisms on all homotopy groups.
\end{defn}

If $k = 1$ this agrees with the usual definition.  Moreover, $\bbU_k$ preserves homotopy groups by a cofinality argument.

\begin{defn}
We say $X \in \sD \sS_k$ is an $\Omega^k$-spectrum if the adjoint structure maps $X(n) \to \Omega^k X(n+1)$ are weak equivalences for all $n$.  A morphism $f \colon X \to Y$ in $\sD \sS_k$ is a stable equivalence if the pullback
\begin{equation*}
f^* \colon [Y, E] \to [X, E]
\end{equation*}
is an isomorphism for every $\Omega^k$-spectrum $E$, where $[-,-]$ indicates the hom-sets in the homotopy category of $\sD \sS_k$ with the level model structure.
\end{defn}

If $\sD = \mathscr{N}$ or $\mathscr{I}$, the stable equivalences and $\pi_*$-isomorphisms coincide.  In the case $\sD = \Sigma$, the stable equivalences strictly contain the $\pi_*$-isomorphisms.  Proofs of these claims are essentially identical to those in \cite[\S 7 \& 8]{mmss} for the case $k=1$.  In fact, everything in \S 7-9 holds essentially verbatim after carefully multiplying all the relevant superscripts by $k$.  As a result we obtain the stable model structure on $k$-$\sD$-spectra.

\begin{defn}
A morphism $f \colon X \to Y$ in $\sD \sS_k$ is an acyclic $q$-cofibration if it is a $q$-cofibration and a stable equivalence.  We say $f$ is a $q$-fibration if it has the right lifting property with respect to all acyclic $q$-cofibrations.
\end{defn}

\begin{thm}
The category $\sD \sS_k$ is a compactly generated proper topological model category with respect to the stable equivalences, $q$-fibrations, and $q$-cofibrations.
\end{thm}

\subsection{Comparison}

We are now in position to state and prove our main comparison result.

\begin{thm}
For $k \geq 1$ and $\sD \in \{ \mathcal{N}, \Sigma, \mathcal{I} \}$, the forgetful functor $\bbU_k \colon \sD \sS \to \sD \sS_k$ is the right adjoint of a Quillen equivalence.  For $\sD = \Sigma$ or $\sD = \mathscr{I}$, the equivalence is symmetric monoidal.
\begin{proof}
For $\sD = \mathscr{N}$ or $\sD = \mathscr{I}$, the result follows from Lemmas 10.2 and 10.3 of \cite{mmss}, since in these cases the stable equivalences and $\pi_*$-isomorphisms coincide.  To handle the case of $\sD = \Sigma$, we use the naturality of $i_k$ to conclude that the following diagram of categories and right Quillen functors commutes.
\begin{equation*}
\begin{tikzcd}
\mathscr{IS} \arrow{r}{\bbU_k} \arrow{d}[swap]{\bbU} & \mathscr{IS}_k \arrow{d}{\bbU} \\
\Sigma \sS \arrow{r}[swap]{\bbU_k} & \Sigma \sS_k
\end{tikzcd}
\end{equation*}
The vertical arrows are Quillen equivalences by arguments in \cite[\S 10]{mmss}, so by 2-of-3 the lower arrow is a Quillen equivalence.
\end{proof}
\end{thm}

\subsection{Model Categories of Modules}

From now on we specialize to $\sD = \Sigma$ or $\sD = \mathscr{I}$.  Just as in the case $k=1$, categories of modules in $\sD \sS_k$ have nice model-theoretic properties.  Compare to \cite[Theorem 12.1]{mmss}.

\begin{thm}
Let $A$ be an $\bbS_k$-algebra.
\begin{enumerate}
\item The category $\sD \sS_A$ of left $A$-modules is a compactly generated proper topological model category with weak equialences and q-fibrations created in $\sD \sS_k$.
\item If $A$ is cofibrant in $\sD \sS_k$, then the forgetful functor $\bbU_A \colon \sD \sS_A \to \sD \sS_k$ preserves q-cofibrations.
\item If $A$ is commutative, the symmetric monoidal category $\sD \sS_A$ satisfies the pushout-product and monoid axioms of \cite{SS00}.
\item If $f \colon B \to A$ is a weak equivalence of $\bbS_k$-algebras, then restriction and extension of scalars define a Quillen equivalence between the categories $\sD \sS_B$ and $\sD \sS_A$.
\end{enumerate}
\end{thm}

\bibliographystyle{amsalpha}
\bibliography{AMRBib}

\end{document}